# The Shortest Path Problem for the Distant Graph of the Projective Line Over the Ring of Integers

## Andrzej Matraś & Artur Siemaszko

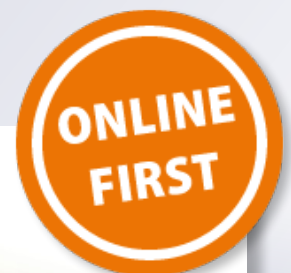
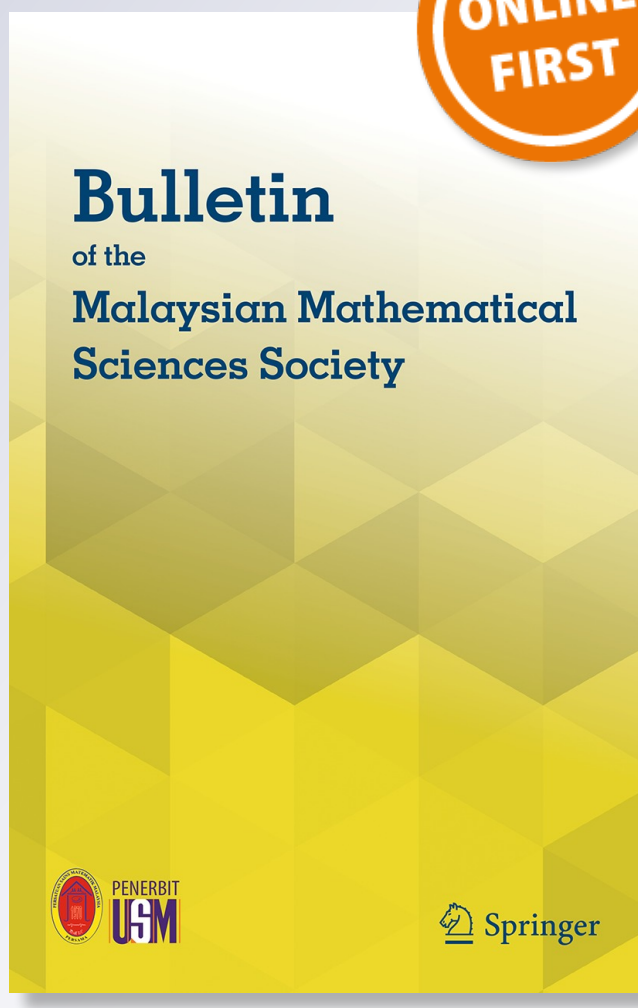

# The Shortest Path Problem for the Distant Graph of the Projective Line Over the Ring of Integers

**Andrzej Matraś**[1] · **Artur Siemaszko**[1]

**Abstract** The distant graph $G = G(\mathbb{P}(Z), \triangle)$ of the projective line over the ring of integers is considered. The shortest path problem in this graph is solved by use of Klein's geometric interpretation of Euclidean continued fractions. In case the minimal path is non-unique, all the possible splitting are described which allows us to give necessary and sufficient conditions for existence of a unique shortest path.

## 1 Introduction

The projective line $\mathbb{P}(R)$ over any ring $R$ with 1 can be defined using admissible pairs [1]. For some rings, for example for commutative ones, it is equivalent to using unimodular pairs.

One of the basics concerning the projective line over ring $\mathbb{P}(R)$ is the distant relation $\triangle$ on the set of points. A pair $(\bar{\mathbf{x}}, \bar{\mathbf{y}}) \in \mathbb{P}(R) \times \mathbb{P}(R)$ is in the distant relation if the matrix $[\mathbf{x}', \mathbf{y}'] \in GL_2(R)$ for any choice of admissible elements $\mathbf{x}' \in \bar{\mathbf{x}}$ and $\mathbf{y}' \in \bar{\mathbf{y}}$. Here, and in the sequel $\bar{\mathbf{x}} := \mathbb{R}\mathbf{x}$ for all admissible $\mathbf{x} \in \mathbb{R}^2$. In that way one defines distant graph $G = G(\mathbb{P}(R), \triangle)$ for any ring $\mathbb{P}(R)$ [4]. The vertices of the distant graph on $\mathbb{P}(R)$ are the points of $\mathbb{P}(R)$, the edges of the graph are the undirected pairs of distant points. The distance function $dist(\bar{\mathbf{x}}, \bar{\mathbf{y}})$ is the minimal number of edges needed to walk from $\bar{\mathbf{x}}$ to $\bar{\mathbf{y}}$. The diameter of the $G$ is the supremum of all distances between its points.

The notion of the distant relation was introduced in the nineties of the previous century (see e.g. [6]). Blunck and Havlicek explored this notion later in [2,3]. They gave for instance the necessary and sufficient condition for the distant graph to be connected.

We consider the projective line over $Z$. In this case, points of $\mathbb{P}(Z)$ are cyclic modules $Z(a, b)$, where $(a, b) \in Z \times Z$ is a unimodular pair. From $Z$ being a $GE_2$ ring follows the connectedness of the distant graph [2], but nothing can be said about its diameter. In this paper, we give a simple direct proof that the distant graph of the projective line over $Z$ has the infinite diameter.

The aim of this paper is to analyse paths in the distant graph $G(\mathbb{P}(Z), \triangle)$ and solve the shortest path problem for it.

The main tool in our considerations is the transition algorithm introduced in [8] (see also [9]) for bases of $Z^2$ in order to investigate some kind of Markov partitions of automorphisms of the two-dimensional torus. This is in fact the reinterpretation of Klein's geometrization of continued fractions (see e.g. [7]).

The distant graph of a ring is defined in terms of "homogenous" ring coordinates. However, our statements follow from affine setup, hence we present in Sect. 2 the transition algorithm on $Z \times Z \subset \mathbb{R}^2$. This is a repetition of [8] but we do it for seek of completeness.

In the next section, we adopt this algorithm to the language of the projective line. Then we make use of it to find the shortest paths connecting given two non-adjacent vertices. It also allows us to provide necessary and sufficient conditions for existence of exactly one shortest path between them.

For any non-adjacent vertices $\bar{\mathbf{x}}, \bar{\mathbf{y}} \in G(\mathbb{P}(Z), \triangle)$ we define so called cone relations on $\mathbb{P}(Z) \setminus \{\bar{\mathbf{x}}, \bar{\mathbf{y}}\}$ with precisely two equivalence classes. Any path entirely contained in one of this classes is called a consistent path. We prove that in each of the equivalence classes of the cone relation there is precisely one shortest path connecting $\bar{\mathbf{x}}$ and $\bar{\mathbf{y}}$ (Theorem 1). All vertices of those two paths together with $\bar{\mathbf{x}}$ and $\bar{\mathbf{y}}$ form the set of vertices of some subgraph, called by us the Klein graph $K(\bar{\mathbf{x}}, \bar{\mathbf{y}})$.

If we abandon the assumption of consistency then the issue becomes more complicated. It turns out that in general there may be lots of shortest paths between two non-adjacent vertices. Nevertheless we prove that there is at least one shortest path which is contained in some induced subgraph of $K(\bar{\mathbf{x}}, \bar{\mathbf{y}})$ called by us a corner graph. This observation allows us to give a recipe to find one of the shortest paths. Consequently, we obtain the formula for the distance $\text{dist}(\bar{\mathbf{x}}, \bar{\mathbf{y}})$ between two vertices (Theorem 2).

Further, the careful insight into the corner graph enables to formulate the necessary and sufficient conditions for the existence of more than one shortest path (Theorem 3).

In [2], Blunck and Havlicek code a path connecting any pair of elements from connected component of the distant graph of a ring $R$ as a product of elementary matrices.

In our modification, this product is a decomposition of an element of $PGL(2, R)$ transforming the first pair of elements of the path onto the last one. In the last section, we get from this factorization the formula for the distance in the distant graph in algebraic language. We hope that this algebraic approach allows us to generalize the results to some wider class of rings, for instance to discretely ordered ones.

## 2 Basic Tools: Klein Sails and Fans of Bases

Formally $\mathbb{P}(Z)$ is a set of cyclic modules $Z(a, b)$ generated by unimodular pairs $(a, b) \in Z \times Z$. To follow all considerations of this paper, it is convenient to lift the projective line to $Z^2$ understood as the subset of the plane $\mathbb{R}^2$.

Every two elements $\bar{\mathbf{x}}, \bar{\mathbf{y}} \in \mathbb{P}(Z)$ determine two lines $k, l \subset \mathbb{R}^2$, respectively. Both lines contain the origin and other integer points, hence, are rational. They divide the plane into four closed quadrants. In each of them consider the boundary of the closed convex hull of all integer points distinct from the origin. This boundary is a infinite polyline. It splits into three parts. Two of those parts are contained in lines $k$ and $l$ and the middle one has precisely one common point with each of those lines (see Fig.1). This part is called a Klein sail. A similar construction may be carried out for irrational lines as well. The difference is that one gets a single polyline having no common points with the lines and asymptotic to them. Klein considered those polylines as geometrization of continued fraction expansions of irrational numbers.

By construction every integer point lying on any of Klein sails form a unimodular pair in $Z \times Z$. We will be interested in the set of those points since every shortest walk between $\bar{\mathbf{x}}$ and $\bar{\mathbf{y}}$ must be contained in the set of $Z$-modules generated by them. The set of integer points on Klein sails may be found by use of the so called transition algorithm developed in [8] (see also [9]) for bases of $Z \times Z$ in order to investigate some kind of Markov partitions of automorphisms of the two-dimensional torus. Fix $\mathbf{x} \in \bar{\mathbf{x}}$ and $\mathbf{y} \in \bar{\mathbf{y}}$ and consider $k$ and $l$ as horizontal and vertical axes positively directed by $\mathbf{x}$ and $\mathbf{y}$. This gives us the enumeration of four quadrants determined by $k$ and $l$. To start the transition algorithm first find, using the Extended Euclidean Algorithm, the bi-sequence $(\mathbf{c}_n)_{n \in Z}$ of all vectors with $\det[\mathbf{x}, \mathbf{c}_n] = 1$ enumerated in such a way that $\mathbf{c}_{n+1} = \mathbf{c}_n - \mathbf{x}$. If $\mathbf{y} \in \{\mathbf{c}_n\}$ then stop the algorithm. Otherwise there is precisely one $n_0 \in Z$ such that $\mathbf{c}_{n_0}$ is in the first quadrant and $\mathbf{c}_{n_0+1}$ is in the second one. Denote $\mathbf{e}_1 := \mathbf{c}_{n_0}$ and $\mathbf{f}_1 := \mathbf{c}_{n_0+1}$. Obviously $\beta_1 = (\mathbf{e}_1, \mathbf{f}_1)$ forms a basis in $Z \times Z$ (see the remark after Lemma 1). Assume the $\beta_k = (\mathbf{e}, \mathbf{f})$ is a basis in $Z \times Z$ with $\mathbf{e}$ in the first quadrant and $\mathbf{f}$ in the second one. If $\mathbf{y} = \mathbf{e} + \mathbf{f}$ then the algorithm stops. Otherwise put $\beta_{k+1} = (\mathbf{e} + \mathbf{f}, \mathbf{f})$ if $\mathbf{e} + \mathbf{f}$ lies in the first quadrant and $\beta_{k+1} = (\mathbf{e}, \mathbf{e} + \mathbf{f})$ lies in the second one.

If the set of bases $\{\beta_k\}$ were infinite then norms of their vectors would tend to infinity while the areas of parallelograms spanned by them would stay equal to 1, hence the angles would tend to 0 or $\pi$. In our construction, they converge to 0 since the angles are decreased step by step. It follows that $\mathbf{y}$ lies in the interior of the parallelogram spanned be vectors of $\beta_k$ for large enough $k$, which is a contradiction (see the remark after Lemma 1). Therefore, the algorithm eventually stops necessarily on $\mathbf{y}$.

In order to translate the above algorithm in terms of the projective line we note that it can be rewritten in the following way. First we may assume that $\mathbf{e}_1 + \mathbf{f}_1$ is in the first quadrant. The other case is completely symmetric. If $\mathbf{e}_1 + \mathbf{f}_1 = \mathbf{y}$ then the algorithm stops. Otherwise we add $\mathbf{f}_1$ to $\mathbf{e}_1$ so many times until the resulting vector stays in the first quadrant. We get the finite sequence of vectors $(\mathbf{e}_{A_0} := \mathbf{e}_1, \mathbf{e}_2, \ldots, \mathbf{e}_{a_1+1} =: \mathbf{e}_{A_1})$. By construction $\mathbf{f}_2 := \mathbf{e}_{A_1} + \mathbf{f}_1$ lies in the second quadrant. Then we keep adding $\mathbf{e}_{A_1}$ to $\mathbf{f}_1$ so many times until the resulting vector stays in the second quadrant. We get the finite sequence of vectors $(\mathbf{f}_{B_0} := \mathbf{f}_1, \mathbf{f}_2, \ldots, \mathbf{f}_{b_1+1} =: \mathbf{f}_{B_1})$. Then we repeat the construction until the summation of vectors is equal to $\mathbf{y}$.

We have constructed two sequences of vectors

$$(\mathbf{e}_1, \mathbf{e}_2, \ldots, \mathbf{e}_{A_1}, \mathbf{e}_{A_1+1}, \ldots, \mathbf{e}_{A_2}, \mathbf{e}_{A_2+1}, \ldots, \mathbf{e}_{A_{r-1}}, \mathbf{e}_{A_{r-1}+1}, \ldots, \mathbf{e}_{A_r})$$

and

$$(\mathbf{f}_1, \mathbf{f}_2, \ldots, \mathbf{f}_{B_1}, \mathbf{f}_{B_1+1}, \ldots, \mathbf{f}_{B_2}, \mathbf{f}_{B_2+1}, \ldots, \mathbf{f}_{B_{l-1}}, \mathbf{f}_{B_{l-1}+1}, \ldots, \mathbf{f}_{B_l}),$$

where $A_k = \sum_{n=0}^{k} a_n$, $k = 0, \ldots, r+1$, and $B_k = \sum_{n=0}^{k} b_n$, $k = 0, \ldots, l+1$. Here $a_0 = a_{r+1} = b_0 = b_{l+1} = 1$. Those two sets together with $\mathbf{x}$ and $\mathbf{y}$ we will project to the projective line and the resulting set of vertices with appropriate edges we will call a *Klein graph* which is the main tool of this paper. The notion of a Klein graph does not depend on choice of representative elements of $\bar{\mathbf{x}}$ and $\bar{\mathbf{y}}$ since opposite Klein sails are symmetric with respect to the origin.

## 3 Basic Tools in Terms of the Projective Line

This section is devoted to presenting the material of the previous one in terms of the projective line. This allows us to conduct all considerations independently on choice of representative elements of vertices in the projective line.

Having two distinct elements $\bar{\mathbf{x}}, \bar{\mathbf{y}} \in \mathbb{P}(Z)$ we define an equivalence relation on $\mathbb{P}(Z) \setminus \{\bar{\mathbf{x}}, \bar{\mathbf{y}}\}$:

$$\bar{\mathbf{u}} \sim_{\{\bar{\mathbf{x}}, \bar{\mathbf{y}}\}} \bar{\mathbf{v}} \quad \text{if} \quad \alpha_{\mathbf{ux}} \alpha_{\mathbf{uy}} \alpha_{\mathbf{vx}} \alpha_{\mathbf{vy}} > 0,$$

for any choice $\mathbf{x} \in \bar{\mathbf{x}}, \mathbf{y} \in \bar{\mathbf{y}}, \mathbf{u} \in \bar{\mathbf{u}}, \mathbf{v} \in \bar{\mathbf{v}}$, where $\alpha_{\mathbf{ux}}, \alpha_{\mathbf{uy}}, \alpha_{\mathbf{vx}}, \alpha_{\mathbf{vy}} \in \mathbb{R}$ such that

$$\mathbf{u} = \alpha_{\mathbf{ux}} \mathbf{x} + \alpha_{\mathbf{uy}} \mathbf{y}, \quad \mathbf{v} = \alpha_{\mathbf{vx}} \mathbf{x} + \alpha_{\mathbf{vy}} \mathbf{y}.$$

Since each of vectors in the inequality appears twice, the relation is well defined. Each such relation is called a *cone relation*.

For two vectors $\mathbf{x}, \mathbf{y} \in Z^2$ denote positive and negative cones in $\mathbb{R}^2$ by

$$C^{\pm}(\mathbf{x}, \mathbf{y}) = \{\alpha \mathbf{x} + \beta \mathbf{y} : \pm \alpha \beta > 0\}.$$

Every cone relation has two equivalence classes being images of

$$C^+(\mathbf{x}, \mathbf{y}) \cap Z^2 \quad \text{and} \quad C^-(\mathbf{x}, \mathbf{y}) \cap Z^2$$

by the canonical projection from $Z^2 \setminus \{(0, 0)\}$ onto the projective line.

The basic property of those relations is the following lemma proved in [8]:

**Lemma 1** (Lemma 2.1, [8]) *If $\bar{\mathbf{x}}, \bar{\mathbf{y}}, \bar{\mathbf{x}}', \bar{\mathbf{y}}' \in \mathbb{P}(Z)$ satisfy $\bar{\mathbf{x}} \triangle \bar{\mathbf{y}}$ and $\bar{\mathbf{x}}' \triangle \bar{\mathbf{y}}'$ then either $\{\bar{\mathbf{x}}, \bar{\mathbf{y}}\} = \{\bar{\mathbf{x}}', \bar{\mathbf{y}}'\}$ or precisely one equivalence class of each cone relations $\sim_{\{\bar{\mathbf{x}},\bar{\mathbf{y}}\}}$ and $\sim_{\{\bar{\mathbf{x}}',\bar{\mathbf{y}}'\}}$ is contained in the equivalence class of another cone relation.*

The above lemma may be expressed in various ways. It is, for example, equivalent to the fact that the interior of the parallelogram spanned by members **e**, **f** of any basis in $Z^2$ contains no vectors from $Z^2$. In other words the only vectors in the closed parallelogram spanned by **e** and **f** are its vertices. By Pick's formula this is equivalent for this parallelogram to have area equal to one. This is in turn equivalent to the equality

$$\det [\mathbf{e}, \mathbf{f}] = \pm 1,$$

where $[\mathbf{e}, \mathbf{f}]$ stands for the matrix with columns equal to **e** and **f** respectively.

We will need also the following obvious observation.

**Lemma 2** *If $\bar{\mathbf{x}}, \bar{\mathbf{y}}, \bar{\mathbf{u}}$ and $\bar{\mathbf{v}}$ are pairwise distinct elements of $\mathbb{P}(Z)$ then*

$$\bar{\mathbf{u}} \sim_{\{\bar{\mathbf{x}},\bar{\mathbf{y}}\}} \bar{\mathbf{v}} \Leftrightarrow \bar{\mathbf{x}} \sim_{\{\bar{\mathbf{u}},\bar{\mathbf{v}}\}} \bar{\mathbf{y}}.$$

In graph theory, a set of vertices is called a *clique* if its elements are pairwise adjacent. From Lemma 1, one can easily deduce the following corollary.

**Corollary 1** 1. *Any maximal clique in $\mathbb{P}(Z)$ has exactly three elements.*
2. *Every pair of adjacent vertices belongs to precisely two maximal cliques.*
3. *If $\{\bar{\mathbf{x}}, \bar{\mathbf{y}}, \bar{\mathbf{u}}\}$ and $\{\bar{\mathbf{x}}, \bar{\mathbf{y}}, \bar{\mathbf{v}}\}$ are two maximal cliques then*

$$\bar{\mathbf{u}} \nsim_{\{\bar{\mathbf{x}},\bar{\mathbf{y}}\}} \bar{\mathbf{v}}.$$

It is not difficult to find maximal cliques for every pair $\bar{\mathbf{x}} \triangle \bar{\mathbf{y}}$. Indeed

$$\{\bar{\mathbf{x}}, \bar{\mathbf{y}}, \overline{\mathbf{x} + \mathbf{y}}\} \quad \text{and} \quad \{\bar{\mathbf{x}}, \bar{\mathbf{y}}, \overline{\mathbf{x} - \mathbf{y}}\}$$

are the only maximal cliques containing $\bar{\mathbf{x}}$ and $\bar{\mathbf{y}}$.

Now we will express the transition algorithm described in the previous section in terms of the projective line. We do this to demonstrate that the notion of Klein graph is really a projective notion as well as to establish some notation used later in the paper. The concept of the Klein graph is a main tool of the next two sections used to find shortest paths connecting an arbitrary pair of vertices of the distant graph.

**Transition Algorithm**

Given two elements $\bar{\mathbf{x}}, \bar{\mathbf{y}} \in \mathbb{P}(Z)$ fix one of them, say $\bar{\mathbf{x}}$.

Using the Extended Euclidean Algorithm find a sequence $(\bar{\mathbf{c}}_n)_{n \in Z}$ of all vertices of the distant graph adjacent to $\bar{\mathbf{x}}$. The sequence $(\bar{\mathbf{c}}_n)_{n \in Z}$ may be ordered in such a way that $\{\bar{\mathbf{c}}_n, \bar{\mathbf{c}}_{n+1}, \bar{\mathbf{x}}\}$ is a clique for every $n \in Z$. If $\bar{\mathbf{y}} \in (\bar{\mathbf{c}}_n)$ then stop.

If not, there is precisely one $n_0 \in Z$ such that $\bar{\mathbf{c}}_{n_0} \nsim_{\{\bar{\mathbf{x}},\bar{\mathbf{y}}\}} \bar{\mathbf{c}}_{n_0+1}$. Denote $\bar{\mathbf{c}}_{n_0}$ and $\bar{\mathbf{c}}_{n_0+1}$ by $\bar{\mathbf{e}}_1$ and $\bar{\mathbf{f}}_1$ arbitrarily.

By 1. and 2. of Corollary 1 there is precisely one $\bar{\mathbf{g}}_2 \neq \bar{\mathbf{x}}$ such that $\{\bar{\mathbf{e}}_1, \bar{\mathbf{f}}_1, \bar{\mathbf{g}}_2\}$ is a clique. Either

$$(\text{TR}) \; \bar{\mathbf{g}}_2 \sim_{\{\bar{\mathbf{x}},\bar{\mathbf{y}}\}} \bar{\mathbf{e}}_1 \quad \text{or} \quad (\text{TL}) \; \bar{\mathbf{g}}_2 \sim_{\{\bar{\mathbf{x}},\bar{\mathbf{y}}\}} \bar{\mathbf{f}}_1.$$

In the sequel of this paper we assume the case (TR). The other case is completely symmetric. It means that considering (TL) instead of (TR) we have to mutually exchange all notations. For instance we exchange symbols $\bar{\mathbf{e}}_j$ and $\bar{\mathbf{f}}_j$. Later on we introduce some numbers $r, a_j, A_j$ associated to vertices $\bar{\mathbf{e}}_j$ and $l, b_j, B_j$ associated to vertices $\bar{\mathbf{f}}_j$. We also must exchange them considering the case (TL) instead of (TR).

Denote $\bar{\mathbf{e}}_2 := \bar{\mathbf{g}}_2$. Among all vertices, excluding $\bar{\mathbf{e}}_1$, adjacent to $\bar{\mathbf{f}}_1$ there are precisely two adjacent to each other that are not in the relation $\sim_{\{\bar{\mathbf{e}}_1,\bar{\mathbf{y}}\}}$. Moreover in the sequence of vertices adjacent to $\bar{\mathbf{f}}_1$ there are finitely many in the relation $\sim_{\{\bar{\mathbf{e}}_1,\bar{\mathbf{y}}\}}$ with $\bar{\mathbf{e}}_2$. We may naturally enumerate them: $(\bar{\mathbf{e}}_2, \ldots, \bar{\mathbf{e}}_{a_1+1})$ in such a way that $\bar{\mathbf{e}}_k \triangle \bar{\mathbf{e}}_{k+1}$. If the next vertex in this sequence and not in the relation $\sim_{\{\bar{\mathbf{e}}_1,\bar{\mathbf{y}}\}}$ with $\bar{\mathbf{e}}_{a_1+1}$ is equal to $\bar{\mathbf{y}}$ then stop.

If not, denote the above vertex by $\bar{\mathbf{f}}_2$ and repeat the same construction with $\bar{\mathbf{e}}_{a_1+1}$ and $\bar{\mathbf{f}}_2$ playing the role of $\bar{\mathbf{f}}_1$ and $\bar{\mathbf{e}}_2$, respectively. We get the sequence $(\bar{\mathbf{f}}_2, \ldots, \bar{\mathbf{f}}_{b_1+1})$. Note that $\bar{\mathbf{f}}_k \triangle \bar{\mathbf{f}}_{k+1}$, $k = 2, \ldots, b_1$, and $\bar{\mathbf{f}}_j \sim_{\{\bar{\mathbf{f}}_1,\bar{\mathbf{y}}\}} \bar{\mathbf{f}}_2$, $j = 2, \ldots, b_1 + 1$.

Then repeat the above construction with $\bar{\mathbf{e}}_{a_1+1}$ and $\bar{\mathbf{f}}_{b_1+1}$ playing the role of $\bar{\mathbf{e}}_1$ and $\bar{\mathbf{f}}_1$, respectively. The transition algorithm is depicted in Fig. 1 for $\mathbf{x} = [2, 3]$ and $\mathbf{y} = [-3, 2]$.

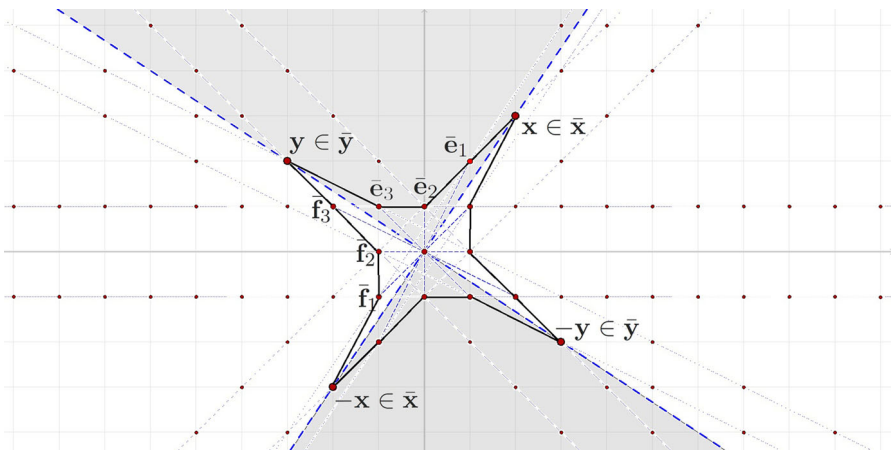

**Fig. 1** Klein sails, Transition algorithm

Since all lines in $\mathbb{R}^2$ determined by members of $\mathbb{P}(Z)$ are rational our algorithm eventually stops after reaching $\bar{\mathbf{y}}$ as we have shown in the previous section. We have constructed two sequences of vertices:

$$(\bar{\mathbf{e}}_1, \bar{\mathbf{e}}_2, \ldots, \bar{\mathbf{e}}_{A_1}, \bar{\mathbf{e}}_{A_1+1}, \ldots, \bar{\mathbf{e}}_{A_2}, \bar{\mathbf{e}}_{A_2+1}, \ldots, \bar{\mathbf{e}}_{A_{r-1}}, \bar{\mathbf{e}}_{A_{r-1}+1}, \ldots, \bar{\mathbf{e}}_{A_r})$$

and

$$(\bar{\mathbf{f}}_1, \bar{\mathbf{f}}_2, \ldots, \bar{\mathbf{f}}_{B_1}, \bar{\mathbf{f}}_{B_1+1}, \ldots, \bar{\mathbf{f}}_{B_2}, \bar{\mathbf{f}}_{B_2+1}, \ldots, \bar{\mathbf{f}}_{B_{l-1}}, \bar{\mathbf{f}}_{B_{l-1}+1}, \ldots, \bar{\mathbf{f}}_{B_l}),$$

where $A_k = \sum_{n=0}^{k} a_n$, $k = 0, \ldots, r+1$, and $B_k = \sum_{n=0}^{k} b_n$, $k = 0, \ldots, l+1$. Here $a_0 = a_{r+1} = b_0 = b_{l+1} = 1$. This is caused by the fact that $\bar{\mathbf{e}}_{A_r}$ and $\bar{\mathbf{f}}_{B_l}$ are both adjacent to $\bar{\mathbf{y}}$ so we also have $\bar{\mathbf{e}}_{A_{r+1}} = \bar{\mathbf{f}}_{B_{l+1}} = \bar{\mathbf{y}}$.

By construction both sequences are contained in different equivalence classes of $\sim_{\{\bar{\mathbf{x}},\bar{\mathbf{y}}\}}$ and two consecutive vertices in each of them are adjacent. Vertices $\bar{\mathbf{e}}_{A_k}$ and $\bar{\mathbf{f}}_{B_k}$, excluded $\bar{\mathbf{y}}$, will be called *corner vertices*. Note that by construction $r - l \in \{0, 1\}$.

Two sequences constructed in the above algorithm may be organized as a sequence of adjacent pairs of vertices.

$$\begin{aligned}\bar{\mathbf{f}}_{B_m} \triangle \bar{\mathbf{e}}_k, \quad & k = A_m, \ldots, A_{m+1}, \quad m = 0, 1, \ldots, r - 1; \\ \bar{\mathbf{e}}_{A_{m+1}} \triangle \bar{\mathbf{f}}_k, \quad & k = B_m, \ldots, B_{m+1}, \quad m = 0, 1, \ldots, l - 1.\end{aligned} \quad (1)$$

To the set of vertices defined above add $\bar{\mathbf{x}}$ and $\bar{\mathbf{y}}$. The resulting set with appropriate edges is an induced subgraph of the distant graph and will be called a *Klein graph* associated to $\bar{\mathbf{x}}$ and $\bar{\mathbf{y}}$. Denote it by $K(\bar{\mathbf{x}}, \bar{\mathbf{y}})$. If we consider only corner vertices in this graph we get an induced subgraph which we will call a *corner graph* and denote by $\widetilde{K}(\bar{\mathbf{x}}, \bar{\mathbf{y}})$.

The Klein graph is depicted in Fig. 2. All vertices of the corner graph are marked in black. The bolded edges are edges of the standard path to be defined in the beginning of Sect. 5.

## 4 Consistent Paths

**Definition 1** A path connecting $\bar{\mathbf{x}}$ and $\bar{\mathbf{y}}$ in the distant graph is called *consistent* if all its elements but $\bar{\mathbf{x}}$ and $\bar{\mathbf{y}}$ are contained in the same equivalence class of $\sim_{\{\bar{\mathbf{x}},\bar{\mathbf{y}}\}}$.

The length of the shortest consistent path between $\bar{\mathbf{x}}$ and $\bar{\mathbf{y}}$ is denoted by $d_c(\bar{\mathbf{x}}, \bar{\mathbf{y}})$.

**Theorem 1** *Let* $\bar{\mathbf{x}}, \bar{\mathbf{y}} \in \mathbb{P}(Z)$.

1. *The union of each equivalence classes of* $\sim_{\{\bar{\mathbf{x}},\bar{\mathbf{y}}\}}$ *and* $\{\bar{\mathbf{x}}, \bar{\mathbf{y}}\}$ *contains a unique consistent path connecting* $\bar{\mathbf{x}}$ *and* $\bar{\mathbf{y}}$. *The lengths of those two paths are equal to*

$$d_a(\bar{\mathbf{x}}, \bar{\mathbf{y}}) = 1 + A_r \quad and \quad d_b(\bar{\mathbf{x}}, \bar{\mathbf{y}}) = 1 + B_l,$$

*where the numbers $A_r$ and $B_l$ are taken from the transition algorithm. In particular the distant graph of $\mathbb{P}(Z)$ is connected.*

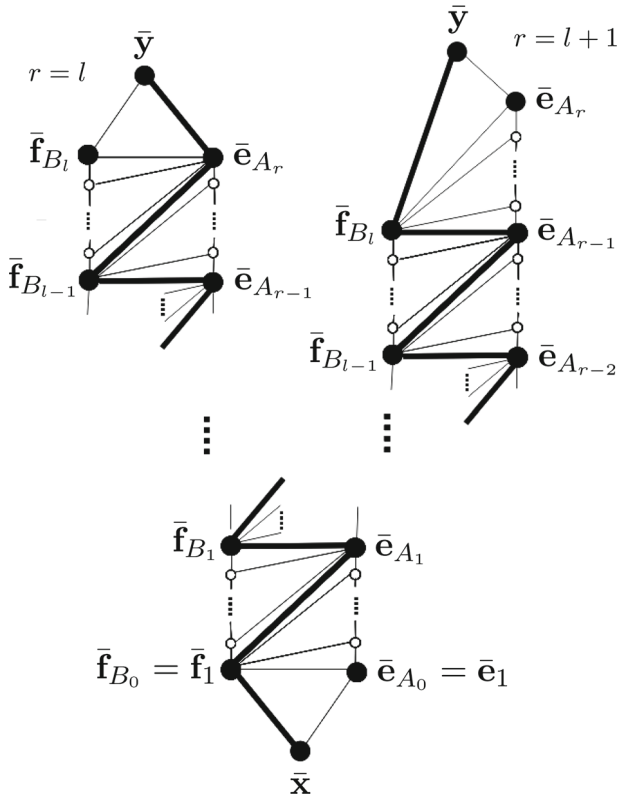

**Fig. 2** Klein graph, corner graph, standard path

2. *There are at most two shortest consistent paths connecting $\bar{\mathbf{x}}$ and $\bar{\mathbf{y}}$ in the distant graph.*
3. *There is a unique shortest consistent path connecting $\bar{\mathbf{x}}$ and $\bar{\mathbf{y}}$ in the distant graph if and only if*

$$d_a(\bar{\mathbf{x}}, \bar{\mathbf{y}}) \neq d_b(\bar{\mathbf{x}}, \bar{\mathbf{y}}).$$

*Proof* Take $\bar{\mathbf{x}}, \bar{\mathbf{y}} \in \mathbb{P}(Z)$. Consider the path

$$p = p(\bar{\mathbf{x}}, \bar{\mathbf{y}}) = \{\bar{\mathbf{x}} = \bar{\mathbf{e}}_0, \bar{\mathbf{e}}_1, \ldots, \bar{\mathbf{e}}_{A_r}, \bar{\mathbf{y}}\},$$

where vertices $\bar{\mathbf{e}}_1, \ldots, \bar{\mathbf{e}}_{A_r}$ are taken from the transition algorithm described in the previous section. It is by construction consistent. Take another consistent path

$$\tilde{p} = \tilde{p}(\bar{\mathbf{x}}, \bar{\mathbf{y}}) = \{\bar{\mathbf{x}} = \bar{\mathbf{z}}_0, \bar{\mathbf{z}}_1 \ldots, \bar{\mathbf{y}}\}$$

contained in the same equivalence class $C$ of $\sim_{\{\bar{\mathbf{x}},\bar{\mathbf{y}}\}}$ as $p$. We will show that $p \subset \tilde{p}$. The set $C\setminus\{\bar{\mathbf{e}}_1, \ldots, \bar{\mathbf{e}}_{A_r}\}$ is a disjoint union of equivalence classes $C_k$ of cone relations $\sim_{\{\bar{e}_k, \bar{e}_{k+1}\}}$, $k = 0, \ldots, A_r$ (here $\bar{e}_0 = \bar{x}$). If $\bar{e}_1 = \bar{z}_1 \in \tilde{p}$ we are done. Otherwise $\bar{z}_1 \in C_0$. It follows from the fact that $\bar{e}_1$ is the last vertex in the sequence of adjacent to $\bar{x}$ vertices contained in $C$ and moreover by Lemma 1, all other vertices in this sequence contained in $C$ belong to $C_0$. We have either $\bar{z}_2 = \bar{e}_1$ or, again by Lemma 1, $\bar{z}_2 \in C_0$. Eventually we get $\bar{e}_1 = \bar{z}_k \in \tilde{p} \cap p$ because we have to go out of $C_0$ in order to get $\bar{y}$. Assuming $\bar{z}_k = \bar{e}_l$ we are able to show, by the same arguments as above, that $\bar{e}_{l+1} \in \tilde{p}$ which shows that $p \subset \tilde{p}$. Therefore $p$ is a unique shortest consistent path between $\bar{x}$ and $\bar{y}$ contained in $C$. When computing the length of $p$ we have to count steps between $\bar{e}_1$ and $\bar{e}_{A_r}$ whose number is $A_r - 1$ and add one step between $\bar{x}$ and $\bar{e}_1$ and another one between $\bar{e}_{A_r}$ and $\bar{y}$. We have finished the proof of 1.

Now 2. and 3. immediately follows from 1. □

*Example 1* Let $\bar{x} = \overline{[1, 0]}$ and $\bar{y} = \overline{[37, 158]}$. Then

$\bar{e}_1 = \overline{[1, 1]}$, $\bar{e}_2 = \overline{[1, 2]}$, $\bar{e}_3 = \overline{[1, 3]}$, $\bar{e}_4 = \overline{[1, 4]}$, $\bar{e}_5 = \overline{[4, 17]}$, $\bar{e}_6 = \overline{[15, 64]}$, $\bar{e}_7 = \overline{[26, 111]}$;
$\bar{f}_1 = \overline{[0, 1]}$, $\bar{f}_2 = \overline{[1, 5]}$, $\bar{f}_3 = \overline{[2, 9]}$, $\bar{f}_4 = \overline{[3, 13]}$, $\bar{f}_5 = \overline{[7, 30]}$, $\bar{f}_6 = \overline{[11, 47]}$.

Therefore $d_a(\bar{x}, \bar{y}) = 8$ and $d_b(\bar{x}, \bar{y}) = d_c(\bar{x}, \bar{y}) = 7$. We see that

$$p(\bar{x}, \bar{y}) = \{\bar{x}, \bar{f}_1, \ldots, \bar{f}_6, \bar{y}\}$$

is the unique shortest consistent path between $\bar{x}$ and $\bar{y}$. Notice also that $\bar{e}_1, \bar{e}_4, \bar{e}_5, \bar{e}_7$; $\bar{f}_1, \bar{f}_4, \bar{f}_6$ are all corner vertices.

*Example 2* Let $\bar{x} = \overline{[1, 0]}$ and $\bar{y} = \overline{[26, 111]}$. Then

$\bar{e}_1 = \overline{[1, 1]}$, $\bar{e}_2 = \overline{[1, 2]}$, $\bar{e}_3 = \overline{[1, 3]}$, $\bar{e}_4 = \overline{[1, 4]}$, $\bar{e}_5 = \overline{[4, 17]}$, $\bar{e}_6 = \overline{[15, 64]}$;
$\bar{f}_1 = \overline{[0, 1]}$, $\bar{f}_2 = \overline{[1, 5]}$, $\bar{f}_3 = \overline{[2, 9]}$, $\bar{f}_4 = \overline{[3, 13]}$, $\bar{f}_5 = \overline{[7, 30]}$, $\bar{f}_6 = \overline{[11, 47]}$.

Therefore $d_a(\bar{x}, \bar{y}) = d_b(\bar{x}, \bar{y}) = d_c(\bar{x}, \bar{y}) = 7$. We see that

$$\{\bar{x}, \bar{e}_1, \ldots, \bar{e}_6, \bar{y}\} \quad \text{and} \quad \{\bar{x}, \bar{f}_1, \ldots, \bar{f}_6, \bar{y}\}$$

are the only two shortest consistent paths between $\bar{x}$ and $\bar{y}$.

**Definition 2** A Hamiltonian cycle $p(\bar{x})$ containing $\bar{x} \in \mathbb{P}(Z)$ is called consistent if there is $\bar{y} \in \mathbb{P}(Z)$ such that $p(\bar{x})$ is a union of two consistent paths connecting $\bar{x}$ and $\bar{y}$.

Formally to run the transition algorithm we have to fix both vertices $\bar{x}$ and $\bar{y}$. But actually we may fix only $\bar{x}$ and arbitrary finite sequences of positive integers $(a_n)$ and $(b_m)$ whose lengths differ by one. Then we may run the transition algorithm starting from the longer sequence. After using all numbers $a_n$ and $b_m$ we get two last corner vertices $\bar{e}_{A_r}$ and $\bar{f}_{B_l}$. Then let $\bar{y}$ be such that it forms a clique with these corner vertices and $\bar{e}_{A_r} \sim_{\{\bar{x}, \bar{y}\}} \bar{f}_{B_l}$. The above consideration gives the following.

**Corollary 2** *For every $\bar{\mathbf{x}} \in \mathbb{P}(Z)$ and every positive integers $d > 2$, $d_a$ and $d_b$ with $d = d_a + d_b$ there exists a Hamiltonian consistent cycle containing $\bar{\mathbf{x}}$ whose length is equal to $d$. Moreover there is $\bar{\mathbf{y}} \in \mathbb{P}(Z)$ such that*

$$d_a = d_a(\bar{\mathbf{x}}, \bar{\mathbf{y}}) \quad \text{and} \quad d_b = d_b(\bar{\mathbf{x}}, \bar{\mathbf{y}}).$$

Note that one can eliminate $\bar{\mathbf{y}}$ from the cycle and get again a Hamiltonian one. In this way one is able to construct decreasing sequence of Hamiltonian cycles containing $\bar{\mathbf{x}}$ that lengths decrease by one. Of course the shortest such cycle has length 3 and forms a clique. In the forthcoming paper it will be shown how to construct all such sequences of Hamiltonian cycles.

## 5 Non-consistent Paths

While, as has been shown in the previous section, there are at most two shortest consistent paths between two vertices in the distant graph of the projective line over $Z$, it turns out that in general you can find many shortest paths if we abandon the assumption of consistency.

In order to compute the length of a shortest path between two vertices $\bar{\mathbf{x}}, \bar{\mathbf{y}} \in \mathbb{P}(Z)$ first run the transition algorithm described in Sect. 3 assuming the case (TR). Then it is convenient to consider the following path which exhaust all vertices of the corner graph:

$$(\bar{\mathbf{x}}, \bar{\mathbf{f}}_1, \bar{\mathbf{e}}_{A_1}, \bar{\mathbf{f}}_{B_1}, \ldots, \bar{\mathbf{f}}_{B_{r-2}}, \bar{\mathbf{e}}_{A_{r-1}}, \bar{\mathbf{f}}_{B_{r-1}}, \bar{\mathbf{e}}_{A_r}, \bar{\mathbf{y}})$$

or

$$(\bar{\mathbf{x}}, \bar{\mathbf{f}}_1, \bar{\mathbf{e}}_{A_1}, \bar{\mathbf{f}}_{B_1}, \ldots, \bar{\mathbf{f}}_{B_{r-2}}, \bar{\mathbf{e}}_{A_{r-1}}, \bar{\mathbf{f}}_{B_{r-1}}, \bar{\mathbf{y}}),$$

depending on whether $r = l$ or $r = l + 1$. Denote this path by $p_s(\bar{\mathbf{x}}, \bar{\mathbf{y}})$ and call it a *standard path* between $\bar{\mathbf{x}}$ and $\bar{\mathbf{y}}$ (see Fig. 2). Let us agree that the vertices $\bar{\mathbf{e}}_{A_k}$ lie on the right side of the corner graph and the vertices of $\bar{\mathbf{f}}_{B_k}$ lie on its left side. The direct insight into the graph $\widetilde{K}(\bar{\mathbf{x}}, \bar{\mathbf{y}})$ gives the following result.

**Lemma 3** *The length of the standard path between non-adjacent $\bar{\mathbf{x}}$ and $\bar{\mathbf{y}}$ is equal to $r + l + 1$ provided $l > 0$. If $l = 0$, then it is equal to 2.*

The key observation of this section is the fact that if want to make a shortest walk from $\bar{\mathbf{x}}$ to $\bar{\mathbf{y}}$ then we have to use only vertices of $K(\bar{\mathbf{x}}, \bar{\mathbf{y}})$ but we may also use only vertices of $\widetilde{K}(\bar{\mathbf{x}}, \bar{\mathbf{y}})$.

**Lemma 4** *Every shortest path between $\bar{\mathbf{x}}$ and $\bar{\mathbf{y}}$ is a subgraph of $K(\bar{\mathbf{x}}, \bar{\mathbf{y}})$. Moreover among all the shortest paths at least one is a subgraph of $\widetilde{K}(\bar{\mathbf{x}}, \bar{\mathbf{y}})$.*

*Proof* We intend to show that having a path $\tilde{p}$ from $\bar{\mathbf{x}}$ to $\bar{\mathbf{y}}$ one can find a path from $\bar{\mathbf{x}}$ to $\bar{\mathbf{y}}$ contained in $K(\bar{\mathbf{x}}, \bar{\mathbf{y}})$ which is not longer than $\tilde{p}$.

Let then $\bar{z} \in \tilde{p} \setminus K(\bar{x}, \bar{y})$. We have either $\bar{z} \sim_{\{\bar{x},\bar{y}\}} \bar{e}_1$ or $\bar{z} \sim_{\{\bar{x},\bar{y}\}} \bar{f}_1$. Assuming the former case we can find precisely one $k$ with $\bar{z} \in C_k$, where $C_k$ are taken from the proof of Theorem 1. By Lemma 1 all vertices adjacent to $\bar{z}$ belong to $C_k \cup \{\bar{e}_k, \bar{e}_{k+1}\}$. Therefore the predecessor as well as the successor of $\bar{z}$ in $\tilde{p}$ are contained in the same set, hence $\{\bar{e}_k, \bar{e}_{k+1}\} \subset \tilde{p}$. Since $\bar{e}_k \triangle \bar{e}_{k+1}$, we may shorten $\tilde{p}$ rejecting all vertices from $\tilde{p} \cap C_k$. We argue analogically in the latter case.

Assume now that $\tilde{p}$ is one of the shortest paths between $\bar{x}$ and $\bar{y}$. We already know that $\tilde{p} \subset K(\bar{x}, \bar{y})$. Let $1 \leq m \leq r$ and $A_m < k < A_{m+1}$ be such that $\bar{e}_k \in \tilde{p}$. In general we have always $\bar{e}_{A_{m+1}} \in \tilde{p}$ since we can not omit $\bar{e}_{A_{m+1}} \in \tilde{p}$ walking from $\bar{e}_k$ to $\bar{y}$. Walking from $\bar{x}$ to $\bar{e}_k$ we must go through $\bar{f}_{B_m}$ or $\bar{e}_{A_m}$ so $\bar{f}_{B_m} \in \tilde{p}$ or $\bar{e}_{A_m} \in \tilde{p}$. Since $\bar{f}_{B_m} \triangle \bar{e}_{A_{m+1}}$ and $\tilde{p}$ is one of the shortest paths, $\bar{f}_{B_m} \notin \tilde{p}$ and $\bar{e}_{A_{m+1}} \in \tilde{p}$. We know that $\bar{e}_{A_m} \triangle \bar{f}_{B_m} \triangle \bar{e}_{A_{m+1}}$ is a path in $\widetilde{K}(\bar{x}, \bar{y})$ of the length 2, hence there are no other vertices between $\bar{e}_{A_m}$ and $\bar{e}_{A_{m+1}}$ but $\bar{e}_k$. We may replace $\bar{e}_{A_m} \triangle \bar{e}_k \triangle \bar{e}_{A_{m+1}}$ by $\bar{e}_{A_m} \triangle \bar{f}_{B_m} \triangle \bar{e}_{A_{m+1}} \subset \widetilde{K}(\bar{x}, \bar{y})$.

We analogically consider $\bar{f}_k \in \tilde{p}$ instead of $\bar{e}_k$ which ends the proof. □

Now let us ponder the question when the standard path is one of the shortest between $\bar{x}$ and $\bar{y}$. Let

$$(a_1, b_1, \ldots) = (c_1, c_2, \ldots, c_{r+l}).$$

The first obvious case is $c_2 = 0$. Indeed, if $c_2 = b_1 = 0$ then $l = 0$ and $\bar{x} \triangle \bar{f}_1 \triangle \bar{y}$ is the standard path of the length 2. If additionally $c_1 = a_1 = 0$ then $r = 0$ and $\bar{x} \triangle \bar{e}_1 \triangle \bar{y}$ is the other shortest path. If $c_2 = 0$ and $c_1 \neq 0$ then $r = 1$ and the standard path is the only shortest one.

The other case is if all numbers $c_k$, $k = 1, \ldots, r+l$, are distinct from 0 (which is equivalent to $c_2 \neq 0$) and all but $c_1$ or $c_{r+l}$ are greater than 1. In such a situation no pairs $(\bar{e}_{A_k}, \bar{e}_{A_{k+1}})$, $k = 1, \ldots, l-1$, $(\bar{f}_{B_k}, \bar{f}_{B_{k+1}})$, $k = 0, \ldots, r-2$, are adjacent. It means that the shortest walk from $\bar{e}_{A_k}$ to $\bar{e}_{A_{k+1}}$ is $\bar{e}_{A_k} \triangle \bar{f}_{B_k} \triangle \bar{e}_{A_{k+1}}$. Similarly the shortest walk between $\bar{f}_{B_k}$ and $\bar{f}_{B_{k+1}}$ is $\bar{f}_{B_k} \triangle \bar{e}_{A_{k+1}} \triangle \bar{f}_{B_{k+1}}$. We have shown the following.

**Proposition 1** *The standard path is one of the shortest ones iff*

- *$c_2 = 0$ or*
- *$c_2 > 0$ and $c_k > 1$ for all $k \in \{2, \ldots, r+l-1\}$.*

Note that even $c_1 = 1$ or $c_{r+l} = 1$ the standard path remains one of the shortest ones. Indeed, if $c_1 = a_1 = 1$ then $\bar{e}_1 \triangle \bar{e}_2$, hence there are two shortest paths from $\bar{x}$ to $\bar{e}_2$, namely $\bar{x} \triangle \bar{f}_1 \triangle \bar{e}_2$ and $\bar{x} \triangle \bar{e}_1 \triangle \bar{e}_2$. Similarly if $c_{r+l-2} = b_{l-1} = 1$ then $\bar{f}_{B_{l-1}} \triangle \bar{f}_{B_l} \triangle \bar{y}$ and $\bar{f}_{B_{l-1}} \triangle \bar{e}_{A_r} \triangle \bar{y}$ are two shortest walks from $\bar{f}_{B_{l-1}}$ to $\bar{y}$ and if $c_{r+l-2} = a_{r-1} = 1$ then $\bar{e}_{A_{r-1}} \triangle \bar{e}_{A_r} \triangle \bar{y}$ and $\bar{e}_{A_{r-1}} \triangle \bar{f}_{B_l} \triangle \bar{y}$ are two shortest walks from $\bar{e}_{A_{r-1}}$ to $\bar{y}$.

One can ask what is a necessary and sufficient condition for the Klein graph $K(\bar{x}, \bar{y})$ so that the standard path is the only shortest one. The above considerations immediately convince us that first of all we must have $c_1 > 1$ and $c_{r+l} > 1$. We must also have $c_k > 1$ for $k = 2, \ldots, r+l-1$ to assure that the standard path is one of the shortest ones. Assume now that $c_k = 2$ for some $k \in \{2, \ldots, r+l-1\}$. If $c_k = a_j$ then $\bar{\mathbf{e}}_{A_{j-1}} \triangle \bar{\mathbf{f}}_{B_{j-1}} \triangle \bar{\mathbf{e}}_{A_j}$ (the standard one) and $\bar{\mathbf{e}}_{A_{j-1}} \triangle \bar{\mathbf{e}}_{A_{j-1}+1} \triangle \bar{\mathbf{e}}_{A_j}$ are two shortest path between $\bar{\mathbf{e}}_{A_{j-1}}$ and $\bar{\mathbf{e}}_{A_j}$. Similarly if $c_k = b_j$ then $\bar{\mathbf{f}}_{B_{j-1}} \triangle \bar{\mathbf{e}}_{A_j} \triangle \bar{\mathbf{f}}_{B_j}$ (the standard one) and $\bar{\mathbf{f}}_{B_{j-1}} \triangle \bar{\mathbf{f}}_{B_{j-1}+1} \triangle \bar{\mathbf{f}}_{B_j}$ are two shortest path between $\bar{\mathbf{f}}_{B_{j-1}}$ and $\bar{\mathbf{f}}_{B_j}$. We have just shown that the necessary and sufficient condition is the following.

**Proposition 2** *The standard path is the only shortest one iff*

- $c_1 > 0$, $c_2 = 0$ *or*
- $c_1 > 1$, $c_{r+l} > 1$ *and* $c_k > 2$ *for all* $k \in \{2, \ldots, r+l-1\}$.

Let us now consider a case when the standard path is not the shortest one. We already know that this is equivalent to

$$c_2 = b_1 \neq 0 \quad \text{and} \quad c_k = 1 \quad \text{for some } k \in \{2, \ldots, r+l-1\}.$$

On the other hand, by Lemma 4, to compute the length of the shortest path between $\bar{\mathbf{x}}$ and $\bar{\mathbf{y}}$ it is enough to find the shortest one contained in the corner graph associated with those two vertices.

The considerations similar to those above allow us to give a recipe for finding one of the shortest paths in $\widetilde{K}(\bar{\mathbf{x}}, \bar{\mathbf{y}})$. First assume the following convention: walking along the standard path we say that we meet the number $a_{k+1}$ or $b_{k+1}$ if we are in the vertex $\bar{\mathbf{e}}_{A_k}$ or $\bar{\mathbf{f}}_{B_k}$, respectively. Then the recipe is the following: start from $\bar{\mathbf{x}}$ and walk along the standard path until you meet the number one. Then walk along the consecutive vertices on the same side of the corner graph. Once you meet the number different from one you come back to the standard path. Perform this procedure till you reach $\bar{\mathbf{e}}_{A_r}$ or $\bar{\mathbf{f}}_{B_l}$. Then make one step more to end at $\bar{\mathbf{y}}$. In the sequel the path we have just described is called a *standard shortest path*.

To compute the distance between $\bar{\mathbf{x}}$ and $\bar{\mathbf{y}}$ it is left to determine by how many steps we shortened the standard path. For the first observe that if $r = l$ then we never meet $a_1$ nor $b_l$ and if $r = l+1$ then we never meet $a_1$ nor $a_r$. Every other appearance of the number one in the "proper" place results in shortening of the standard path by precisely one step. It is possible to shorten the standard path by one step in precisely two cases:

1. $b_k = 1$;    2. $b_k > 1$ and $a_k = 1$.

The first case always works. To decide whether 2. works or not we must check if it satisfies some additional condition. Put

$$s_k = \begin{cases} \max\{j : b_{k-1} = \ldots = b_{k-j} = 1\} : b_{k-1} = 1, \\ 0 : \hspace{10em} b_{k-1} > 1; \end{cases} \quad (2)$$

and

$$t_k = \begin{cases} \max\{j : a_{k-1} = \ldots = a_{k-j} = 1\} : a_{k-1} = 1, \\ 0 : \hspace{10em} a_{k-1} > 1. \end{cases} \quad (3)$$

Then 2. results in shortening iff $s_k \leq t_k$.

Now we are in position to determine by how many steps we shorten the standard path. To end this define

$$\tilde{a}_k = \begin{cases} a_k : & k = 2, \ldots, l; \\ \max(a_k, 2) : & k = 1 \end{cases} \quad (4)$$

and

$$\tilde{b}_k = \begin{cases} b_k : & k = 1, \ldots, l-1; \\ \max(b_k, 2) : & k = l, \end{cases} \quad (5)$$

if $r = l$. If $r = l + 1$ the former condition remains unchanged and the latter one becomes

$$\tilde{b}_k = b_k \quad \text{for} \quad k = 1, \ldots, l.$$

Now remove from $\{1, \ldots, l\}$ those $k$'s that satisfy 2. with $s_k > t_k$ and denote the received set by $D = D_{\bar{x},\bar{y}}$. By the above consideration we see that we shorten the standard path by

$$\sum_{k \in D} \left\lfloor \frac{1}{2} \left( \left\lfloor \frac{1}{\tilde{a}_k} \right\rfloor + \left\lfloor \frac{1}{\tilde{b}_k} \right\rfloor \right) + \frac{1}{2} \right\rfloor$$

steps. The symbol $\lfloor \cdot \rfloor$ stands for the floor function.

**Theorem 2** *Assume that $\bar{x}, \bar{y} \in \mathbb{P}(Z)$ are not adjacent. Then the length of the shortest path connecting $\bar{x}$ and $\bar{y}$ is equal to*

$$\mathrm{dist}(\bar{x}, \bar{y}) = \begin{cases} r + l + 1 - \sum_{k \in D_{\bar{x},\bar{y}}} \left\lfloor \frac{1}{2} \left( \left\lfloor \frac{1}{\tilde{a}_k} \right\rfloor + \left\lfloor \frac{1}{\tilde{b}_k} \right\rfloor \right) + \frac{1}{2} \right\rfloor : l > 0 \\ 2 : & l = 0. \end{cases}$$

*In particular, $G(\mathbb{P}(Z), \Delta)$ is of infinite diameter.*

## 6 Uniqueness of the Shortest Paths

In this section we formulate necessary and sufficient conditions for the existence of more than one shortest path between two vertices. In the previous section we described the construction of one of the shortest paths between two non-adjacent vertices. This path may split, but need not to, into more shortest paths because of appearance of 1 or 2 in the sequence $(c_1, \ldots, c_{r+l})$. Two next lemmas explain the situation.

First we deal with ones. We say the two blocks of 1's *overlap* if there are natural numbers $0 \leq l_1 \leq l_2 \leq l + 1$ and $0 \leq r_1 \leq r_2 \leq r + 1$ satisfying

(o1) $l_1 < r_1 \leq l_2 + 1 \leq r_2$ with $l_2 < l + 1$ iff $r = l + 1$

or

(o2) $l_1 < r_1 \leq l_2 = r_2 = r + 1 = l + 1$

or

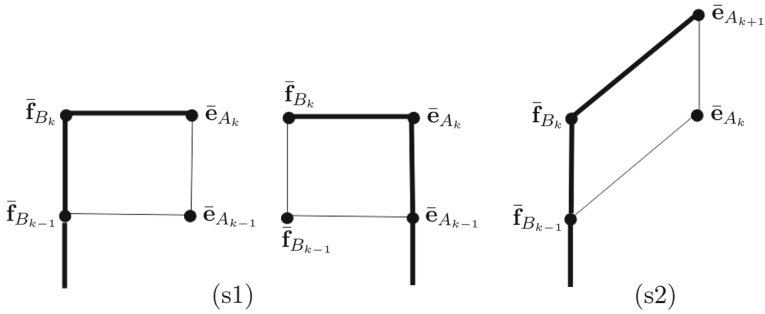

**Fig. 3** Splitting by 1

(o3) $0 = r_1 = l_1 \leq l_2 < r_2$

or

(o4) $r_1 \leq l_1 \leq r_2 \leq l_2$, where $l_1 \neq 0$ with an additional condition $r_2 < r + 1$ iff $r = l$

or

(o5) $r_1 \leq l_1 \leq l_2 < r_2 = r + 1$.

such that

$$b_{l_1} = b_{l_1+1} = \ldots = b_{l_2} = a_{r_1} = a_{r_1+1} = \ldots = a_{r_2} = 1.$$

**Lemma 5** *Appearance of 1 in the sequence* $(c_1, \ldots, c_{r+l})$ *results in splitting the shortest path iff two blocks of 1's overlap or* $r = l$ *and* $c_i = 1$ *for all i's.*

*Proof* First note that 1 in $(c_1, \ldots, c_{r+l})$ causes splitting of any path $p$ contained in the Klein graph iff it appears in one of the following configurations depicted in Fig. 3:

(s1) $b_k = a_k = 1, 1 \leq k \leq l$;
(s2) $b_k = 1, a_{k+1} = 1, \bar{\mathbf{f}}_k \in p$.

In particular to split the shortest path 1 has to appear in both sides of the Klein graph. Then fix some block $B$ of 1's in the left side of the Klein graph: $b_{l_1} = \cdots = b_{l_2} = 1$, where $0 < l_1 \leq l_2 < l$. Then imagine that we "shift" some block $A$ of 1's in the right side: $a_{r_1} = \ldots = a_{r_2} = 1$, where $0 < r_1 \leq r_2 < r$. We want to assure that meeting of those two blocs results in splitting the shortest path.

Assume that $|A| = r_2 - r_1 + 1 \leq |B| = l_1 - l_2 + 1$. Then in order to gain the configuration (s1) or (s2) that causes splitting we must have

$$r_1 \leq l_2 + 1 \leq r_2 \quad \text{or} \quad r_1 \leq l_1 \leq r_2.$$

Indeed, if neither of the above conditions holds then there are three possibilities. Two of them are $l_2 + 1 < r_1$ or $l_1 > r_2$. In this cases neither (s1) nor (s2) holds. The last case is $r_1 > l_1$ and $l_2 \geq r_2$. This forces $|A| < |B|$ and appearance of (s1) and (s2). Nevertheless this does not cause splitting because the standard shortest path has as a part the path $\bar{\mathbf{f}}_{B_{l_1-1}} \triangle \bar{\mathbf{f}}_{B_{l_1}} \triangle \cdots \triangle \bar{\mathbf{f}}_{B_{l_2}}$. Walking along this part we do not "jump" onto the right side of the Klein graph although we meet 1's in configuration (s1) and (s2). Therefore we are not able to gain splitting. Together with $|A| \leq |B|$ the former condition gives (o1) and the latter one gives (o4).

Similarly, if $|A| \geq |B|$ then we must have

$$l_1 < r_1 \leq l_2 + 1 \quad \text{or} \quad l_1 \leq r_2 \leq l_2.$$

Together with $|A| \geq |B|$ the former condition gives (o1) and the latter one gives (o4).

Until now we have assumed that our blocks neither start nor end at $\bar{\mathbf{x}}$ and $\bar{\mathbf{y}}$. Now we assume opposite.

Assume first $r_1 = l_1 = 0$. Then in order to have splitting we must have $r_2 > l_2$ which gives (o3). If $l_1 = 0 \neq r_1$ then $l_1 < r_1$ and necessarily $r_1 \leq l_2 + 1 \leq r_2$. This gives (o1). If $r_1 = 0 \neq l_1$ then $l_1 > r_1$ and necessarily $l_1 \leq r_2 \leq l_2$ which is (o4).

Next let $r = l$ and $r_2 = l_2 = r$. Then in order to have splitting we must have $l_1 < r_1$ which is (o2). If $l_2 = l = r \neq r_2$ then $r_2 < r$ which is equivalent to $r_2 < r + 1$. Moreover we must have $r_1 \leq l_1 \leq r_2$. This gives the additional condition in (o4). If $l_2 \neq l = r = r_2$ then $l_2 + 1 \leq r_2$ and we must have $l_1 < r_1 \leq l_2 + 1$ which is (o1).

Now assume that $r = l + 1$. First consider $r_2 < r$. Then $r_2 < r = l + 1 = l_2$ and in order to guarantee splitting we must have $r_1 \leq l_1 \leq r_2$. This gives (o4). If $r_2 = r + 1$ then to obtain splitting the standard shortest path has to end in the right side of the Klein graph. Therefore if $r_1 \leq l_1$ then $l_2 \leq l + 1 < r_2 = l + 2$ and we get (o5). On the other hand if $r_1 > l_1$ then we must have $r_1 \leq l_2 + 1$ and we have to "jump" into the right side of the Klein graph which forces $l_2 < l + 1$. This is the additional condition in (o1).

There is left one case. Namely, $c_1 = \cdots = c_{r+l} = 1$. In this case it obvious that one has splitting iff $r = l$.

We have considered all possible configurations of two blocks of 1's that result in splitting. Since we get all conditions (o1)–(o5) and each of them obviously forces splitting, we have finished the proof. □

Another reason of splitting is appearance of 2 in this sequence. By direct insight into the graph $K(\bar{\mathbf{x}}, \bar{\mathbf{y}})$ we get the following result.

**Lemma 6** *Appearance of 2 in the sequence* $(c_1, \ldots, c_{r+l})$ *results in splitting the standard shortest path iff one of the following conditions hold*

(st1) *there exists* $k \in \{2, \ldots, l\}$ *such that* $b_k = 2$, $a_k > 1$, $a_{k+1} > 1$ *or* $b_1 = 2$ *and* $a_2 > 1$.
(st2) *there exists* $k \in \{2, \ldots, l\}$ *such that* $a_k = 2$, $b_k > 1$ *and* $s_k \leq t_k$,

*where the numbers $s_k$ and $t_k$ in the above are defined in* (2) *and* (3).

*Proof* Note that if $b_k = 2$ then necessary and sufficient condition to split the standard shortest path is that it contains vertices $\bar{\mathbf{f}}_{B_{k-1}}, \bar{\mathbf{f}}_{B_k}$ and $\bar{\mathbf{e}}_{A_k}$ (see the left side of Fig. 4). The standard shortest path contains $\bar{\mathbf{f}}_{B_{k-1}}$ iff $k > 1$ and $a_k > 1$ or $k = 1$ even if $a_1 = 1$. It necessarily contains $\bar{\mathbf{e}}_{A_k}$ since $b_k = 2 > 1$. The standard shortest path contains $\bar{\mathbf{f}}_{B_k}$ iff $a_{k+1} > 1$.

**Fig. 4** Splitting by 2

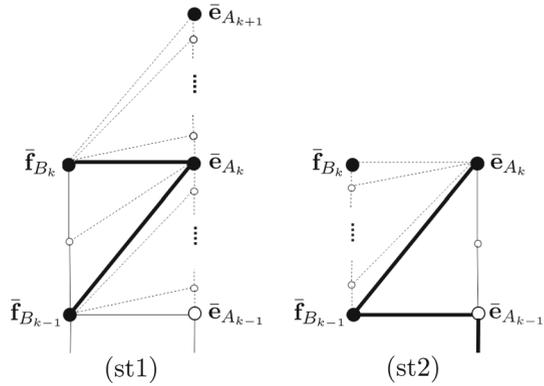

Now assume that $a_k = 2$. Then necessary and sufficient condition to split the standard shortest path is appearance of vertices $\bar{\mathbf{e}}_{A_{k-1}}, \bar{\mathbf{f}}_{B_{k-1}}$ and $\bar{\mathbf{e}}_{A_k}$ (see the right side of Fig. 4). The vertex $\bar{\mathbf{e}}_{A_0} = \bar{\mathbf{e}}_1$ is never contained in the standard shortest path so $k$ must be greater than 1. The vertex $\bar{\mathbf{e}}_{A_{k-1}}$ appears in it iff $s_k \leq t_k$. This path necessarily contains $\bar{\mathbf{f}}_{B_{k-1}}$ since $a_k = 2 > 1$. Then we see that $\bar{\mathbf{e}}_{A_k}$ appears in the standard shortest path iff $b_k > 1$.

Observe that in the both cases $k$ must be less then $r + 1$. This gives range of $k$ in the thesis of the lemma. □

All discussion of this section and direct insight into the Klein graph convince us that the following statement is true.

**Theorem 3** *Assume that $\bar{\mathbf{x}}, \bar{\mathbf{y}} \in \mathbb{P}(Z)$ are not adjacent. There is more than one shortest path between $\bar{\mathbf{x}}$ and $\bar{\mathbf{y}}$ if and only if*

– *two blocks of 1's overlap or*
– *$r = l$ and all numbers $a_k$ and $b_k$ are equal to one or*
– *one of the conditions of Lemma 6 holds.*

## 7 Algebraic Representation of Paths in the Klein Graph

In [2] there is described a one-to-one correspondence between products of matrices of the form

$$E(a) = \begin{pmatrix} a & 1 \\ -1 & 0 \end{pmatrix} \tag{6}$$

and paths in $G = G(\mathbb{P}(Z), \triangle)$. The length of a path is equal to the number of matrices in an associated product. This together with some relations in $GE_2(Z)$ noted in [5] gives the algebraic method to find the shortest walks between arbitrary vertices of the distant graph.

For any two non distant vertices $\bar{\mathbf{x}}$ and $\bar{\mathbf{y}}$ of the distant graph we express the method of finding the shortest path developed in Sect. 5 in terms of matrices form $GE_2(Z)$. Fix two representative elements $\mathbf{x} \in \bar{\mathbf{x}}$ and $\mathbf{y} \in \bar{\mathbf{y}}$. In order to use the Blunck–Havlicek algorithm we establish as a standard basis in $Z^2$ the basis $(\mathbf{x}, \mathbf{f}_1)$, where $\mathbf{f}_1$ is defined in Sect. 2 for $\mathbf{x}$ and $\mathbf{y}$. Let $\mathbf{y}$ represents as $[p, q]$ in the above basis. Then by construction $q/p \geq 1$. The (TR) case of Sect. 3 requires $\mathbf{y}$ to have the slope not less then 2. If $1 \leq q/p < 2$ then we change the basis to $(-\mathbf{x}, \mathbf{e}_1)$. In this basis $\mathbf{y}$ represents as $[q - p, q]$, hence it has the slope not less than 2.

The standard path from $\bar{\mathbf{x}}$ to $\bar{\mathbf{y}}$ is associated by the Blunck–Havlicek algorithm with the matrix

$$E\left((-1)^{n+1}(d_n)\right) \cdot \ldots \cdot E(d_1) \cdot E(-d_0) \cdot E(0), \tag{7}$$

where $q/p = [d_0; d_1, \ldots, d_n]$ is a continued fraction expansion with $d_n > 1$.

The number of matrices in the product in (7) may be reduced using the following two formulas

$$E(a)E(\pm 1)E(b) = \pm E(a \mp 1)(b \mp 1). \tag{8}$$

There are two more formulas reducing a number of matrices in any product of matrices of the form (6), namely

$$E(a)E(\pm 2)^2 E(b) = \mp E(a \mp 1)E(\mp 3)E(\mp 1). \tag{9}$$

Nevertheless, it is not difficult to see that after all possible reductions of (7) using (8) the expression from the left hand side of (9) does not appear in the resulting product.

Any change of order of reductions may effect the resulting representation of $A$ as a product of matrices of the form (6). In spite of this all representations that cannot be reduced any more have the same number of matrices since they all represent the shortest paths from $\bar{\mathbf{x}}$ to $\bar{\mathbf{y}}$.

We say that a matrix $A$ *represents a path* from $\bar{\mathbf{x}}$ to $\bar{\mathbf{y}}$ if in the basis described above for any choice $\mathbf{x} \in \bar{\mathbf{x}}$ to $\mathbf{y} \in \bar{\mathbf{y}}$ the product given by the Blunck–Havlicek algorithm is equal to $A$.

By results of Sect. 5 the matrix given by (7) represents at least one shortest path from $\bar{\mathbf{x}}$ to $\bar{\mathbf{y}}$. There may be other matrices representing such a path. The situations is described in the below proposition whose thesis directly follow from analysis of the Klein graph.

**Proposition 3** *There are at most two matrices representing the shortest paths from $\bar{\mathbf{x}}$ to $\bar{\mathbf{y}}$, among them there is always a matrix $A$ given by* (7).

*All the shortest paths from $\bar{\mathbf{x}}$ to $\bar{\mathbf{y}}$ are represented by precisely one matrix iff*

1. *either $d_n > 2$ or*
2. *$d_n - 1 = d_{n-1} = \ldots = d_{n-k} = 1$ and $d_{n-k-1} > 1$ for some odd $k \geq 1$.*

*If there are two matrices representing the shortest paths from $\bar{\mathbf{x}}$ to $\bar{\mathbf{y}}$ then the other one is equal to*

$$-\begin{pmatrix} 1 & 0 \\ (-1)^n & 1 \end{pmatrix} \cdot A.$$

**Acknowledgments** We would like to express our gratitude to professor Hans Havlicek for attracting our attention to the problem and fruitful discussion on the subject. We also thank to Mariusz Bodzioch, a Ph.D. student of our faculty, for help in preparing the figures. The second named author was partially supported by the NCN (National Science Centre, Poland) Grant 2011/03/B/ST1/04427.